\documentclass[twoside,11pt,reqno]{amsart}
\usepackage{amsmath,amsthm,amssymb,amstext,amsfonts,amscd}
\usepackage{graphicx}\usepackage{multirow}\usepackage[cp1250]{inputenc}\usepackage{ifthen}
\newtheorem{theorem}{Theorem}[section]
\newtheorem{definition}{Definition}[section]
\newtheorem{example}{Example}[section]\newtheorem{lemma}{Lemma}[section]
\newtheorem{remark}{Remark}[section]
\newenvironment{Proof}[1][Proof]{\noindent\textbf{#1.} }{\ \rule{0.0em}{0.0em}}
\numberwithin{equation}{section}

\begin{document}
\title[\textbf{Growth of composite entire functions....}]{\textbf{Measure of
relative $(p,q)$-th order based on a growth of composite entire
functions }}
\author[S. Kanas, S. K. Datta, T. Biswas, G. K. Mondal]{S. Kanas$%
^{\ast }$, S. K. Datta, T. Biswas and G. K. Mondal}
\address{S. Kanas : Department of Mathematical Analysis, Faculty of Mathematics
and Natural Sciences, University of Rzeszow, Poland}
\email{skanas@ur.edu.pl}
\address{S. K. Datta : Department of Mathematics, University of Kalyani,
P.O.- Kalyani, Dist-Nadia, PIN-\ 741235, West Bengal, India}
\email{sanjib\_kr\_datta@yahoo.co.in}
\address{T. Biswas : Rajbari, Rabindrapalli, R. N. Tagore Road, P.O.-
Krishnagar, Dist-Nadia, PIN-\ 741101, West Bengal, India}
\email{tanmaybiswas\_math@rediffmail.com}
\address{G. K. Mondal : Dhulauri Rabindra Vidyaniketan (H.S.), Vill +P.O.-
Dhulauri , P.S.- Domkal, Dist-Murshidabad , PIN-\ 742308, West Bengal, India}
\email{golok.mondal13@rediffmail.com}
\thanks{$^{\ast }$ Corresponding author: E-Mail: skanas@ur.edu.pl}



\begin{abstract}
We deduce some growth properties of composite entire functions in
the light of their relative $(p,q)$\ th order by extending some
results of J. Tu, Z. X. Chen and X. M. Zheng \cite{14}.
\end{abstract}

\maketitle

\section{\textbf{Background, fundamental definitions and notations}}

Let $f$ be an entire function defined on  a set of all complex
numbers $\mathbb{C}$.
The maximum modulus function $M_{f}$ or $M_f(r)$ of $f(z)=\overset{\infty }{\underset{n=0%
}{\sum }}a_{n}z^{n}$ on $\left\vert z\right\vert =r$ is defined as $M_{f}=%
\underset{\left\vert z\right\vert =r}{\max }\left\vert f\left(
z\right) \right\vert $. If $f$ is non-constant entire, then its
maximum modulus function $M_{f}\left( r\right) $ is strictly
increasing and continuous, and therefore there exists its inverse
function $M_{f}^{-1}:\left( \left\vert f\left( 0\right) \right\vert
,\infty \right) \rightarrow \left( 0,\infty \right) $ with
$\underset{s\rightarrow \infty }{\lim }M_{f}^{-1}\left( s\right)
=\infty .$ Moreover, for given any two entire functions $f$ and $g$
the ratio $\frac{M_{f}\left( r\right) }{M_{g}\left( r\right) }$, as
$r\rightarrow \infty$, is called \textit{the growth of $f$ with
respect to $g$} in terms of their maximum moduli. Our notations are
standard within the theory of Nevanlinna's value distribution of
entire functions, and therefore we do not explain those in detail as
available in \cite{15}. In the sequel the following two notations
are used:
\begin{eqnarray*}
\log ^{[k]}x &=&\log \left( \log ^{[k-1]}x\right) \text{ for
}k=1,2,3, \cdots, \\
\log ^{[0]}x &=&x,
\end{eqnarray*}%
and%
\begin{eqnarray*}
\exp ^{[k]}x &=&\exp \left( \exp ^{[k-1]}x\right) \text{ for
}k=1,2,3, \cdots, \\
\exp ^{[0]}x &=&x.
\end{eqnarray*}

Let us recall that Juneja, Kapoor and Bajpai \cite{10} defined the $%
(p,q)$\emph{-th order} and $(p,q)$\emph{-th lower order},
respectively, of an entire function $f$ as follows:

\begin{equation*}
\rho _{f}\left( p,q\right) =\text{ }\underset{r\rightarrow \infty
}{\lim \sup }\frac{\log ^{[p]}M_{f}(r)}{\log ^{\left[ q\right]
}r},\quad \text{and} \quad \lambda _{f}\left( p,q\right) =\text{
}\underset{r\rightarrow \infty }{\lim \inf }\frac{\log
^{[p]}M_{f}(r)}{\log ^{\left[ q\right] }r},
\end{equation*}%
where $p,q$ are positive integers with $p\geq q$.

In this connection we just recall the following definition:

\begin{definition}
\label{d1}\cite{10} An entire function $f$ is said to have
index-pair $\left( p,q\right) $, $p\geq q\geq 1$ if $b<\rho
_{f}\left( p,q\right) <\infty $ and $\rho _{f}\left( p-1,q-1\right)
$ is not a nonzero finite number, where $b=1$ if $p=q$, and $b=0$ if
$p>q.$ Moreover, if $0<\rho _{f}\left( p,q\right) <\infty ,$ then
\begin{equation*}\left\{\begin{array}{lcl}
\rho _{f}\left( p-n,q\right) =\infty& \textit{for}&n<p,\\
~\rho_{f}\left( p,q-n\right) =0&\textit{for}&n<q,\\
\rho _{f}\left( p+n,q+n\right)
=1&\textit{for}&n=1,2,....~.\end{array}\right.
\end{equation*}
\end{definition}

Similarly for $0<\lambda _{f}\left( p,q\right) <\infty$, one can
easily verify that%
\begin{equation*}\left\{\begin{array}{lcl}
\lambda _{f}\left( p-n,q\right) =\infty & \textit{for}&n<p,\\
\lambda_{f}\left( p,q-n\right) =0& \textit{for}&n<q,\\
\lambda _{f}\left( p+n,q+n\right) =1& \textit{for}&n=1,2,....~.
\end{array}\right.\end{equation*}

The definition of $(p,q)$\emph{-th order} ($(p,q)$\emph{%
-th lower order}, respectively), as initiated by Juneja, Kapoor and
Bajpai \cite{10}, extends the notion of \emph{generalized order
}$\rho _{f}^{\left[ l\right] }$ (\emph{generalized lower order
}$\lambda _{f}^{\left[ l\right]}$, resp.) of an entire function $f$
introduced by Sato \cite{12} for each
integer $l\geq 2$, as these correspond to the particular case $\rho _{f}^{%
\left[ l\right] }=\rho _{f}\left(l,1\right) $ ($\lambda _{f}^{%
\left[ l\right] }=\lambda _{f}\left( l,1\right) $, resp. )$.$ If
$p=2$ and $q=1$, then we write $\rho _{f}\left( 2,1\right) =\rho
_{f}$ ($\lambda _{f}\left( 2,1\right) =\lambda _{f}$, resp.) which
is known as \emph{order}
(\emph{lower order}, resp.) of an entire function $f$. The \emph{order%
} (\emph{lower order}, resp.) of an entire function $f$ is classical
in complex analysis and is generally used in computational purpose
which is defined in terms of the growth of $f$ with respect to the
function $\exp z$
function as:%
\begin{equation*}
\rho _{f}=\underset{r\rightarrow \infty }{\lim \sup }\frac{\log \log
M_{f}\left( r\right) }{\log \log M_{\exp z}\left( r\right) }=\underset{%
r\rightarrow \infty }{\lim \sup }\frac{\log \log M_{f}\left( r\right) }{\log
r}
\end{equation*}%
\begin{equation*}
\left(\lambda _{f}=\underset{r\rightarrow \infty }{\lim \inf
}\frac{\log \log M_{f}\left( r\right) }{\log \log M_{\exp z}\left(
r\right) }=\underset{r\rightarrow \infty }{\lim \inf }\frac{\log
\log M_{f}\left( r\right) }{\log r}, \text{resp.}\right) .
\end{equation*}

Bernal \cite{1, 2} introduced the relative order between two entire
functions to avoid comparing growth just with $\exp z$ which is as
follows:

\begin{definition}
\label{d2}\cite{1, 2} The relative order of $f$ with respect to $g$,
denoted as $\rho _{g}\left( f\right)$, is defined
by:%
\begin{eqnarray*}
\rho _{g}\left( f\right) &=&\inf \left\{ \mu >0:M_{f}\left( r\right)
<M_{g}\left( r^{\mu }\right) \text{ for all }r>r_{0}\left( \mu \right)
>0\right\} \\
&=&\underset{r\rightarrow \infty }{\lim \sup }\frac{\log
M_{g}^{-1}M_{f}\left( r\right) }{\log r}.
\end{eqnarray*}
\end{definition}

This definition coincides with the classical one if $g=\exp z$
\cite{13}.  Similarly, one can define the relative lower order of
$f$ with respect to $g$ denoted by $\lambda _{g}\left( f\right) $ as
\begin{equation*}
\lambda _{g}\left( f\right) =\underset{r\rightarrow \infty }{\lim \inf }%
\frac{\log M_{g}^{-1}M_{f}\left( r\right) }{\log r}.
\end{equation*}

Lahiri and Banerjee \cite{11} gave a more generalized concept of
relative order in the following way:

\begin{definition}
\label{d3}\cite{11} If $k\geq 1$ is a positive integer, then the $k$-th
generalized relative order of $f$ with respect to $g$, denoted by $\rho
_{f}^{k}\left( g\right) $ is defined by%
\begin{eqnarray*}
\rho _{g}^{k}\left( f\right) &=&\inf \left\{ \mu >0:M_{f}\left( r\right)
<M_{g}\left( \exp ^{\left[ k-1\right] }r^{\mu }\right) \text{ for all }%
r>r_{0}\left( \mu \right) >0\right\} \\
&=&\underset{r\rightarrow \infty }{\lim \sup }\frac{\log ^{\left[ k\right]
}M_{g}^{-1}M_{f}\left( r\right) }{\log r}.
\end{eqnarray*}%
Clearly, $\rho _{g}^{1}\left( f\right) =\rho _{g}\left( f\right) $ and $\rho
_{\exp }^{1}\left( f\right) =\rho _{f}$.
\end{definition}

The following definition of relative $(p,q)$ th order of an entire
function in the light of index-pair is due to Sanchez Ruiz et. al.
\cite{y}:

\begin{definition}
\label{d4}\cite{y} Let $f$ and $g$ be any two entire functions with index-pairs $%
\left( m,q\right)$ (and $\left( m,p\right) $ resp.) where $p,q,m$
are positive integers such that $m\geq \max (p,q).$ Then the
relative $\left( p,q\right) $-th order of $f$ with respect to $g$ is
defined as
\begin{equation*}
\rho _{g}^{\left( p,q\right) }\left( f\right) =\underset{r\rightarrow \infty
}{\lim \sup }\frac{\log ^{\left[ p\right] }M_{g}^{-1}M_{f}\left( r\right) }{%
\log ^{\left[ q\right] }r}.
\end{equation*}%
The relative $\left( p,q\right) $-th lower order of $f$ with respect to $g$
is defined by:
\begin{equation*}
\lambda _{g}^{\left( p,q\right) }\left( f\right) =\underset{r\rightarrow
\infty }{\lim \inf }\frac{\log ^{\left[ p\right] }M_{g}^{-1}M_{f}\left(
r\right) }{\log ^{\left[ q\right] }r}.
\end{equation*}
\end{definition}

The previous definitions are easily generated from above as
particular cases, e.g. if $f$ and $g$ have got index-pair $\left(
m,1\right) $ and $\left( m,k\right)$, resp., then Definition
\ref{d4} reduces to Definition \ref{d3}. If the entire functions $f$
and $g$ have the same index-pair $\left( p,1\right)$, where $p$ is
any positive integer, we get
the definition of relative order introduced by Bernal \cite{1}, and if $%
g=\exp ^{\left[ m-1\right] }z,$ then $\rho _{g}\left( f\right) =\rho _{f}^{%
\left[ m\right] }$ and $\rho _{g}^{\left( p,q\right) }\left( f\right) =\rho
_{f}\left( m,q\right) .$ And, if $f$ is an entire function with index-pair $%
\left( 2,1\right) $ and $g=\exp z$, then Definition \ref{d4} becomes the
classical one given in \cite{13}.\medskip

In order to calculate the growth rates of entire functions, the
notions of use of\ the growth indicators such as \emph{order} and \emph{%
lower order} are classical in complex analysis and during the past
decades, several researchers have already been continuing their
studies in the area of comparative growth properties of composite
entire functions in different directions using the classical growth
indicators. But at that time, the concepts of \emph{relative orders}
and \emph{relative lower orders} of entire functions as well as
their technical advantages of not comparing with the growths of
$\exp z$ are not at all known to the researchers of this area.
Therefore the studies of the growths of composite entire functions
in the light of their relative orders and relative lower orders are
the prime concern of this paper. In fact, some light has already
been thrown on such type of works by Datta et. al. in \cite{4, 5, 6}
and \cite{7}. Taking into account all these above, we discuss in
this paper some growth
properties of composite entire functions in the light of their \emph{%
relative }$(p,q)$\emph{\ th order} and \emph{relative
}$(p,q)$\emph{\ th lower order}, after improving some results of J.
Tu, Z. X. Chen and X. M. Zheng \cite{14}\emph{.}

\section{\textbf{Some examples}}

In this section we present some examples of entire functions in
connection with definitions given in the previous section.

\begin{example}[Order of $\exp $]
Given any natural number $m$, the exponential function $f(z)=\exp z^{m}$ has
got $M_{f}\left( r\right) =\exp r^{m}$. Therefore $\frac{\log ^{\left[ 2%
\right] }M_{f}\left( r\right) }{\log r}$ is constantly equal to $m$ and
consequently,
\begin{equation*}
\rho _{f}=\lambda _{f}=m.
\end{equation*}
\end{example}

\begin{example}[Generalized order]
Given any natural numbers $l,m$, the function $f(z)=\exp ^{[l]}z^{m}$ has
got $M_{f}\left( r\right) =\exp ^{[l]}r^{m}$. Therefore $\frac{\log ^{\left[
k\right] }M_{f}\left( r\right) }{\log r}$ is constant for each natural $%
k\geq 2$, thereby following that
\begin{equation*}
\rho _{f}^{\left[ l+1\right] }=\lambda _{f}^{\left[ l+1\right] }=m,
\end{equation*}%
but $\rho _{f}^{\left[ k\right] }=\lambda _{f}^{\left[ k\right] }=+\infty $
for $2\leq k\leq l$, and $\rho _{f}^{\left[ k\right] }=\lambda _{f}^{\left[ k%
\right] }=0$ for $k>l+1.$
\end{example}

\begin{example}[Index-pair]
Given any four positive integers $k,n,p,q$ with $p\geq q$, the function $%
f(z)=\exp ^{[k]}z^{n}$ generates a constant quotient $\frac{\log ^{\left[ p%
\right] }M_{f}\left( r\right) }{\log ^{\left[ q\right] }r}$, and
clearly
\begin{equation*}
\rho _{f}\left( p,q\right) =\lambda _{f}\left( p,q\right) =n\text{ for }%
(p,q)=(k+1,1),
\end{equation*}%
but
\begin{equation*}
\rho _{f}\left( p,q\right) =\lambda _{f}\left( p,q\right) =\left\{
\begin{array}{lcl}
1&\text{ for }&(p,q)=(k+h,h)=1, h\in \mathbb{N}, \\
\infty &\text{for}&p\leq q+1, \\
0&\text{for}&p\geq q+1.%
\end{array}\right.\end{equation*}
Thus $f$ is a regular function with growth $(k+1,1)$.
\end{example}

\begin{example}[Relative (p, q)-th order between functions]
Suppose $f(z)$ $=$ $\exp ^{k}\left\{ z^{n}\right\} $ and $g(z)=\exp ^{\left[ k%
\right] }\left\{ z^{m}\right\} $ with $k,m,n$ any three positive integers.
Then $f$ and $g$ are regular functions with $(k+1,1)$-growth with
\begin{equation*}
\rho _{f}\left( k+1,1\right) =n,\quad \rho _{g}\left( k+1,1\right) =m.
\end{equation*}%
In order to find out their $(1,1)$ relative order we evaluate that
\begin{equation*}
\frac{\log M_{g}^{-1}M_{f}\left( r\right) }{\log r}=\frac{\log \frac{1}{m}%
\left\{ \log ^{\left[ k\right] }\left( \exp ^{\left[ k\right] }r^{n}\right)
\right\} ^{\frac{1}{m}}}{\log r}
\end{equation*}%
which happens to be constant. By taking limits, we easily get
\begin{equation*}
\rho _{g}^{\left( 1,1\right) }\left( f\right) =\lambda _{g}^{\left(
1,1\right) }\left( f\right) =\frac{n}{m}.
\end{equation*}
\end{example}

\section{\textbf{Growth of composite entire functions}}

First of all, we recall one related known property which will be
needed in order to prove our results, as we see in the following
lemma.

\begin{lemma}
\label{l1}\cite{3} If $f$ and $g$ are two entire functions, then for
all sufficiently large values of $r$
\begin{equation*}
M_{f}\left( \frac{1}{8}M_{g}\left( \frac{r}{2}\right) -\left\vert
g\left( 0\right) \right\vert \right) \leq M_{f\circ g}(r)\leq
M_{f}\left( M_{g}\left( r\right) \right).
\end{equation*}
\end{lemma}

Now we present the main results concerning the growth of the
composite entire functions $f$ and $g$.

\begin{theorem}
\label{t1} Let $f$ and $g$ be any two entire functions with index-pairs $%
\left( p,q\right) $ and $\left( m,n\right)$, resp., where $p,q,m,n$
are all positive integers such that $p\geq q$ and $m\geq n.$
Then\newline $\noindent $\textbf{(i)} the index-pair of $f\circ g$
is $\left( p,n\right) $ when $q=m$ and either $\lambda _{f}\left(
p,q\right) >0$ or $\lambda _{g}\left( m,n\right) >0.$ Also
\begin{eqnarray*}
\left( a\right) ~\lambda _{f}\left( p,q\right) \rho _{g}\left(
m,n\right) &\leq &\rho _{f\circ g}\left( p,n\right) \leq \rho
_{f}\left( p,q\right) \rho _{g}\left( m,n\right) \text{ if }\lambda
_{f}\left(
p,q\right) >0,\text{ and} \\
\left( b\right) ~\lambda _{f}\left( p,q\right) \rho _{g}\left(
m,n\right) &\leq &\rho _{f\circ g}\left( p,n\right) \leq \rho
_{f}\left( p,q\right) \rho _{g}\left( m,n\right) \text{ if }\lambda
_{g}\left( m,n\right) >0;
\end{eqnarray*}%
$\noindent $\textbf{(ii)} the index-pair of $f\circ g$ is $\left(
p,q+n-m\right) $ when $q>m$, and either $\lambda _{f}\left(
p,q\right) >0$ or
$\lambda _{g}\left( m,n\right) >0.$ Also%

\noindent$\left( a\right) ~\lambda _{f}\left( p,q\right) \leq \rho
_{f\circ g}\left( p,q+n-m\right) \leq \rho _{f}\left( p,q\right)
\text{ if }\lambda _{f}\left( p,q\right) >0,\text{ ~and}$

\noindent$ \left( b\right) ~\rho _{f\circ g}\left( p,q+n-m\right)
=\rho _{f}\left( p,q\right) \text{ if }\lambda _{g}\left( m,n\right)
>0;$

\noindent\textbf{(iii)} the index-pair of $f\circ g$ is $\left(
p+m-q,n\right) $ when $q<m$, and either $\lambda _{f}\left(
p,q\right) >0$ or $\lambda _{g}\left( m,n\right) >0.$ Also

\noindent $\left( a\right) ~\rho _{f\circ g}\left( p+m-q,n\right)
=\rho _{g}\left( m,n\right) \text{ if }\lambda _{f}\left( p,q\right)
>0$, and

\noindent $\left( b\right) ~\lambda _{g}\left( m,n\right) \leq \rho
_{f\circ g}\left( p+m-q,n\right) \leq \rho _{g}\left( m,n\right)
\text{ if }\lambda _{g}\left( m,n\right) >0.$

\begin{Proof}\rm
In view of the first part of Lemma \ref{l1}, it follows for all sufficiently
large values of $r$ that%
\begin{equation}
\log ^{\left[ p\right] }M_{f\circ g}\left( r\right) \geq \left(
\lambda _{f}\left( p,q\right) -\varepsilon \right) \log ^{\left[
q\right] }M_{g}\left( \frac{r}{2}\right) +O(1),  \label{2}
\end{equation}%
and also for a sequence of values of $r$ tending to infinity we have
\begin{equation}
\log ^{\left[ p\right] }M_{f\circ g}\left( r\right) \geq \left( \rho
_{f}\left( p,q\right) -\varepsilon \right) \log ^{\left[ q\right]
}M_{g}\left( \frac{r}{2}\right) +O(1).  \label{3}
\end{equation}%
Similarly, in view of the second part of Lemma \ref{l1}, for all
sufficiently large values of $r$ we obtain
\begin{equation}
\log ^{\left[ p\right] }M_{f\circ g}\left( r\right) \leq \left( \rho
_{f}\left( p,q\right) +\varepsilon \right) \log ^{\left[ q\right]
}M_{g}\left( r\right). \label{1}
\end{equation}%
Now, the following two cases may arise:
\newline \textbf{Case I.}
$q=m$.
\newline From  \eqref{1} for all sufficiently large values
of $r$, we have
\begin{eqnarray}
\log ^{\left[ p\right] }M_{f\circ g}\left( r\right)  &\leq &\left(
\rho _{f}\left( p,q\right) +\varepsilon \right) \left( \rho
_{g}\left( m,n\right)
+\varepsilon \right)   \notag \\
\textrm{so\ that}~\underset{r\rightarrow \infty }{\lim }\frac{\log
^{\left[ p\right] }M_{f\circ g}\left( r\right) }{\log ^{\left[
n\right] }r} &\leq &\rho _{f}\left( p,q\right) \rho _{g}\left(
m,n\right).\label{4}
\end{eqnarray}%
Also from \eqref{2}, for a sequence of values of $r$ tending to
infinity, we obtain
\begin{eqnarray}
\log ^{\left[ p\right] }M_{f\circ g}\left( r\right)  &\geq &\left( \lambda
_{f}\left( p,q\right) -\varepsilon \right) \left( \rho _{g}\left( m,n\right)
-\varepsilon \right) \log ^{\left[ n\right] }r+O(1),\ \textrm{hence}\notag \\
\underset{r\rightarrow \infty }{\lim \sup }\frac{\log ^{\left[ p\right]
}M_{f\circ g}\left( r\right) }{\log ^{\left[ n\right] }r} &\geq &\lambda
_{f}\left( p,q\right) \rho _{g}\left( m,n\right) \text{~}.  \label{5}
\end{eqnarray}%
Moreover, from \eqref{3} for a sequence of values of $%
r$ tending to infinity, we get
\begin{eqnarray}
\log ^{\left[ p\right] }M_{f\circ g}\left( r\right)  &\geq &\left( \rho
_{f}\left( p,q\right) -\varepsilon \right) \left( \lambda _{g}\left(
m,n\right) -\varepsilon \right) \log ^{\left[ n\right] }r+O(1),\ \textrm{and} \notag \\
\underset{r\rightarrow \infty }{\lim \sup }\frac{\log ^{\left[
p\right] }M_{f\circ g}\left( r\right) }{\log ^{\left[ n\right] }r}
&\geq &\rho _{f}\left( p,q\right) \lambda _{g}\left( m,n\right).
\label{6}
\end{eqnarray}%
Therefore for $\lambda _{f}\left( p,q\right) >0$ and from \eqref{4}
and \eqref{5}, we see  that
\begin{eqnarray}
\lambda _{f}\left( p,q\right) \rho _{g}\left( m,n\right)  &\leq &\underset{%
r\rightarrow \infty }{\lim \sup }\frac{\log ^{\left[ p\right]
}M_{f\circ g}\left( r\right) }{\log ^{\left[ n\right] }r}\leq \rho
_{f}\left(
p,q\right) \rho _{g}\left( m,n\right),   \notag \\
i.e.,~\lambda _{f}\left( p,q\right) \rho _{g}\left( m,n\right) &\leq
&\rho _{f\circ g}\left( p,n\right) \leq \rho _{f}\left( p,q\right)
\rho _{g}\left( m,n\right). \label{7}
\end{eqnarray}%
Likewise, \eqref{4} and \eqref{6} for $\lambda _{g}\left( m,n\right)
>0$ yields
\begin{eqnarray}
\rho _{f}\left( p,q\right) \lambda _{g}\left( m,n\right)  &\leq &\underset{%
r\rightarrow \infty }{\lim \sup }\frac{\log ^{\left[ p\right]
}M_{f\circ g}\left( r\right) }{\log ^{\left[ n\right] }r}\leq \rho
_{f}\left(
p,q\right) \rho _{g}\left( m,n\right)   \notag \\
i.e.,~\rho _{f}\left( p,q\right) \lambda _{g}\left( m,n\right) &\leq
&\rho _{f\circ g}\left( p,n\right) \leq \rho _{f}\left( p,q\right)
\rho _{g}\left( m,n\right). \label{8}
\end{eqnarray}%
Also from \eqref{7} and \eqref{8} one can easily verify that $\rho
_{f\circ g}\left( p-1,n\right) =\infty ,$ $\rho _{f\circ g}\left(
p,n-1\right) $ $=$ $0$ and $\rho _{f\circ g}\left( p+1,n+1\right)
=1$, and therefore we obtain that the index-pair of $f\circ g$ is
$\left( p,n\right) $ when $q=m$, and either $\lambda _{f}\left(
p,q\right)
>0$ or $\lambda _{g}\left( m,n\right) >0$. Thus the first part of the
theorem is established.
\newline \textbf{Case II.}  $q>m$.
\newline
Now, from \eqref{1} for all sufficiently large values of $r$, we
obtain
\begin{eqnarray*}
\log ^{\left[ p\right] }M_{f\circ g}\left( r\right)  &\leq&\left(
\rho
_{f}\left( p,q\right) +\varepsilon \right) \log ^{\left[ q-m\right] }\log ^{%
\left[ m\right] }M_{g}\left( r\right)   \notag \\
i.e.,~\log ^{\left[ p\right] }M_{f\circ g}\left( r\right)  &\leq
&\left( \rho _{f}\left( p,q\right) +\varepsilon \right) \log ^{\left[ q-m%
\right] }\left[ \left( \rho _{g}\left( m,n\right) +\varepsilon \right) \log
^{\left[ n\right] }r\right]   \notag \\
i.e.,~\log ^{\left[ p\right] }M_{f\circ g}\left( r\right)  &\leq
&\left( \rho _{f}\left( p,q\right) +\varepsilon \right) \log ^{\left[ q+n-m%
\right] }r+O(1),\end{eqnarray*} therefore
\begin{equation}\label{9}~\underset{r\rightarrow \infty }{\lim }\frac{\log ^{\left[
p\right] }M_{f\circ g}\left( r\right) }{\log ^{\left[ q+n-m\right]
}r} \leq \rho _{f}\left( p,q\right).
\end{equation}
Also, from \eqref{2} for a sequence of values of $r$ tending to
infinity, we have
\begin{eqnarray*}
\log ^{\left[ p\right] }M_{f\circ g}\left( r\right)  &\geq &\left( \lambda
_{f}\left( p,q\right) -\varepsilon \right) \log ^{\left[ q-m\right] }\left[
\left( \rho _{g}\left( m,n\right) -\varepsilon \right) \log ^{\left[ n\right]
}\left( \frac{r}{2}\right) \right]  \\
&&~\ \ \ \ \ \ \ \ \ \ \ \ \ \ \ \ \ \ \ \ \ \ \ \ \ \ \ \ \ \ \ \ \ \ \ \ \
\ \ \ \ \ \ +O(1)
\end{eqnarray*}
\begin{eqnarray*}
i.e.,~\log ^{\left[ p\right] }M_{f\circ g}\left( r\right)  &\geq
&\left( \lambda _{f}\left( p,q\right) -\varepsilon \right) \log
^{\left[ q-m+n\right] }r+O(1), \notag \end{eqnarray*}hence
\begin{equation}\label{10}\underset{r\rightarrow \infty }{\lim \sup }\frac{\log ^{\left[
p\right] }M_{f\circ g}\left( r\right) }{\log ^{\left[ q+n-m\right]
}r} \geq \lambda _{f}\left( p,q\right).
\end{equation}
Further, for a sequence of values of $r$ tending to infinity,
\eqref{3} yields
\begin{eqnarray*} \log ^{\left[ p\right] }M_{f\circ
g}\left( r\right)  &\geq &\left( \rho _{f}\left( p,q\right)
-\varepsilon \right) \log ^{\left[ q-m\right] }\left[
\left( \lambda _{g}\left( m,n\right) -\varepsilon \right) \log ^{\left[ n%
\right] }\left( \frac{r}{2}\right) \right]  \\
&&~\ \ \ \ \ \ \ \ \ \ \ \ \ \ \ \ \ \ \ \ \ \ \ \ \ \ \ \ \ \ \ \ \ \ \ \ \
\ \ \ \ \ \ +O(1)
\end{eqnarray*}
\begin{eqnarray*}
i.e.,~\log ^{\left[ p\right] }M_{f\circ g}\left( r\right)  &\geq
&\left( \rho _{f}\left( p,q\right) -\varepsilon \right) \log
^{\left[ q+n-m\right] }r+O(1),\end{eqnarray*} so that
\begin{equation}\label{11}\underset{r\rightarrow \infty }{\lim \sup
}\frac{\log ^{\left[ p\right] }M_{f\circ g}\left( r\right) }{\log
^{\left[ q+n-m\right] }r} \geq \rho _{f}\left(
p,q\right).\end{equation} Therefore, from \eqref{9} and \eqref{10}
for $\lambda _{f}\left( p,q\right) >0$, we obtain
\begin{eqnarray}
\lambda _{f}\left( p,q\right)  &\leq &\underset{r\rightarrow \infty }{\lim
\sup }\frac{\log ^{\left[ p\right] }M_{f\circ g}\left( r\right) }{\log ^{%
\left[ q+n-m\right] }r}\leq \rho _{f}\left( p,q\right)   \notag \\
i.e.,~\lambda _{f}\left( p,q\right)  &\leq &\rho _{f\circ g}\left(
p,q+n-m\right) \leq \rho _{f}\left( p,q\right). \label{12}
\end{eqnarray}%
Likewise, for $\lambda _{g}\left( m,n\right) >0$, \eqref{9} and
\eqref{11} follows
\begin{eqnarray}
\rho _{f}\left( p,q\right)  &\leq &\underset{r\rightarrow \infty }{\lim \sup
}\frac{\log ^{\left[ p\right] }M_{f\circ g}\left( r\right) }{\log ^{\left[
q+n-m\right] }r}\leq \rho _{f}\left( p,q\right)   \notag \\
i.e.,~\rho _{f\circ g}\left( p,q+n-m\right)  &=&\rho _{f}\left(
p,q\right). \label{13}
\end{eqnarray}%
Hence, from \eqref{12} and \eqref{13}, one can easily verify that
$\rho _{f\circ g}\left( p-1,q+n-m\right) =\infty ,$ $\rho _{f\circ
g}\left( p,q+n-m-1\right) =0$, and $\rho _{f\circ g}\left(
p+1,q+n-m+1\right) =1$. Therefore we get that the index-pair of
$f\circ g$ is $\left( p,q+n-m\right) $ when $q>m$ and either
$\lambda _{f}\left( p,q\right) >0$ or $\lambda _{g}\left( m,n\right)
>0$, and thus the second part of the theorem follows.\newline
\textbf{Case III.}  $q<m$.
\newline For all sufficiently large values
of $r$ and by \eqref{1} we obtain
\begin{eqnarray*}
\log ^{\left[ p+m-q\right] }M_{f\circ g}\left( r\right)  &\leq &\log ^{%
\left[ m\right] }M_{g}\left( r\right) +O(1)  \notag \\
i.e.,~\log ^{\left[ p+m-q\right] }M_{f\circ g}\left( r\right) &\leq
&\left( \rho _{g}\left( m,n\right) +\varepsilon \right) \log ^{\left[ n%
\right] }r+O(1),\end{eqnarray*} so that
\begin{equation}\underset{r\rightarrow \infty }{\lim }\frac{\log ^{\left[
p+m-q\right] }M_{f\circ g}\left( r\right) }{\log ^{\left[ n\right]
}r} \le \rho _{g}\left( m,n\right). \label{14}
\end{equation}
Also, from \eqref{2}  for a sequence of values of $r$ tending to
infinity, we have
\begin{eqnarray*}
\log ^{\left[ p+m-q\right] }M_{f\circ g}\left( r\right)  &\geq &\log ^{\left[
m\right] }M_{g}\left( \frac{r}{2}\right) +O(1)  \notag \\
i.e.,~\log ^{\left[ p+m-q\right] }M_{f\circ g}\left( r\right)  &\geq
&\left( \rho _{g}\left( m,n\right) -\varepsilon \right) \log
^{\left[ n\right] }r+O(1),\end{eqnarray*}therefore \begin{equation}
\underset{r\rightarrow \infty }{\lim \sup }\frac{\log ^{\left[
p+m-q\right] }M_{f\circ g}\left( r\right) }{\log ^{\left[ n\right]
}r} \geq \rho _{g}\left( m,n\right). \label{15}
\end{equation}
Further, an application of \eqref{3} for a sequence of values of $r$
tending to infinity gives
\begin{eqnarray*}
\log ^{\left[ p+m-q\right] }M_{f\circ g}\left( r\right)  &\geq &\log ^{\left[
m\right] }M_{g}\left( \frac{r}{2}\right) +O(1)  \notag \\
i.e.,~\log ^{\left[ p+m-q\right] }M_{f\circ g}\left( r\right)  &\geq
&\left( \lambda _{g}\left( m,n\right) -\varepsilon \right) \log
^{\left[ n\right] }r+O(1),\end{eqnarray*} and so
\begin{equation}\underset{r\rightarrow \infty }{\lim \sup }\frac{\log ^{\left[
p+m-q\right] }M_{f\circ g}\left( r\right) }{\log ^{\left[ n\right]
}r} \ge \lambda _{g}\left( m,n\right). \label{16}
\end{equation}
Therefore, \eqref{14} and \eqref{15} applied for $\lambda _{f}\left(
p,q\right) >0$ implies
\begin{eqnarray}
\rho _{g}\left( m,n\right)  &\leq &\frac{\log ^{\left[ p+m-q\right]
}M_{f\circ g}\left( r\right) }{\log ^{\left[ n\right] }r}\leq \rho
_{g}\left( m,n\right)   \notag \\
i.e.,~\rho _{f\circ g}\left( p+m-q,n\right)  &=&\rho _{g}\left(
m,n\right). \label{17}
\end{eqnarray}%
Similarly, \eqref{14} and \eqref{16} for $\lambda _{g}\left(
m,n\right) >0$ yields
\begin{eqnarray}
\lambda _{g}\left( m,n\right)  &\leq &\underset{r\rightarrow \infty }{\lim
\sup }\frac{\log ^{\left[ p+m-q\right] }M_{f\circ g}\left( r\right) }{\log ^{%
\left[ n\right] }r}\leq \rho _{g}\left( m,n\right)   \notag \\
i.e.,~\lambda _{g}\left( m,n\right)  &\leq &\rho _{f\circ g}\left(
p+m-q,n\right) \leq \rho _{g}\left( m,n\right). \label{18}
\end{eqnarray}%
An application of  the relation \eqref{17} and \eqref{18} easily
gives that $\rho _{f\circ g}\left( p+m-q-1,n\right) =\infty$,  $\rho
_{f\circ g}\left( p+m-q,n-1\right) =0$ and $\rho _{f\circ g}\left(
p+m-q+1,n+1\right) =1$. Therefore we obtain that the index-pair of $%
f\circ g$ is $\left( p+m-q,n\right) $ when $q<m$ and either $\lambda
_{f}\left( p,q\right) >0$ or $\lambda _{g}\left( m,n\right) >0$, and
thus the third part of the theorem is established.
\end{Proof}{\hfill{$\Box$}{\medskip}}
\end{theorem}

\begin{remark}
\label{r1} Theorem \ref{t1} can be treated as an extension of Theorem $3.1$
and Theorem $3.2$ of Tu, Chen and Zheng \cite{14}.
\end{remark}

\begin{theorem}
\label{t2} Let $f$ and $g$ be any two entire functions with index-pairs $%
\left( p,q\right) $ and $\left( m,n\right)$, resp., where $p,q,m,n$
are all positive integers such that $p\geq q$ and $m\geq n.$ Then

\begin{eqnarray*}
\left( i\right) ~\lambda _{f}\left( p,q\right) \lambda _{g}\left( m,n\right)
&\leq &\lambda _{f\circ g}\left( p,n\right)  \\
&\leq &\min \left\{ \rho _{f}\left( p,q\right) \lambda _{g}\left(
m,n\right),\lambda _{f}\left( p,q\right) \rho _{g}\left( m,n\right)
\right\}
\end{eqnarray*}%
\begin{equation*}
~\ \ \ \ \ \ \ \ \ \ \ \ \ \ \ \ \ \ \ \ \ \ \ \ \ \ \ \ \ \ \ \ \ \
\ \ \ \ \ \text{if }q=m,\text{ }\lambda _{f}\left( p,q\right)
>0\text{ and }\lambda _{g}\left( m,n\right) >0,
\end{equation*}%
\begin{equation*}
\left( ii\right) ~\lambda _{f\circ g}\left( p,q+n-m\right) =\lambda
_{f}\left( p,q\right) \text{ if }q>m,\text{ }\lambda _{f}\left( p,q\right) >0%
\text{ and }\lambda _{g}\left( m,n\right) >0,
\end{equation*}%
and%
\begin{equation*}
\left( iii\right) ~\lambda _{f\circ g}\left( p+m-q,n\right) =\lambda
_{g}\left( m,n\right) \text{ if }q<m\text{, }\lambda _{f}\left( p,q\right) >0%
\text{ and }\lambda _{g}\left( m,n\right) >0.
\end{equation*}
\end{theorem}

Reasoning similarly as in the proof of the Theorem \ref{t1} one can
easily deduce the conclusion of Theorem \ref{t2}, and so its proof
is omitted.

\begin{theorem}
\label{t3} Let $f,g,h$ and $k$ be any four entire functions with
index-pairs $\left( p,q\right) $, $\left( m,n\right)$, $\left(
a,b\right) $ and $\left( c,d\right)$, resp.,  where
$a,b,c,d,p,q,m,n$ are all positive integers such that $a\geq b$,
$c\geq d,$ $p\geq q$ and $m\geq n.$
\newline
$\noindent $\textbf{(i)} If either ($q=m,$ $a=c=p,$ $q\geq n$) or ($q<m,$ $%
c=p,$ $a=p+m-q,$ $q\geq n$) holds and $\lambda _{f}\left( p,q\right) >0$, $%
0<\lambda _{h}^{\left( b,n\right) }\left( f\circ g\right) \leq \rho
_{h}^{\left( b,n\right) }\left( f\circ g\right) <\infty $, $0<\lambda
_{k}^{\left( d,q\right) }\left( f\right) \leq \rho _{k}^{\left( d,q\right)
}\left( f\right) <\infty$, then%
\begin{equation*}
\frac{\lambda _{h}^{\left( b,n\right) }\left( f\circ g\right) }{\rho
_{k}^{\left( d,q\right) }\left( f\right) }\leq \underset{r\rightarrow \infty
}{\lim \inf }\frac{\log ^{\left[ b\right] }M_{h}^{-1}M_{f\circ g}\left(
r\right) }{\log ^{\left[ d\right] }M_{k}^{-1}M_{f}\left( \exp ^{\left[ q-n%
\right] }r\right) }\leq \frac{\lambda _{h}^{\left( b,n\right) }\left( f\circ
g\right) }{\lambda _{k}^{\left( d,q\right) }\left( f\right) }~~\ \ \ \ \ \ \
\ \ \ \
\end{equation*}%
\begin{equation*}
~\ \ \ \ \ \ \ \ \ \ \ \ \ \ \ \ \ \ \ \ \ \ \ \ \ \ \ \ \ \leq \underset{%
r\rightarrow \infty }{\lim \sup }\frac{\log ^{\left[ b\right]
}M_{h}^{-1}M_{f\circ g}\left( r\right) }{\log ^{\left[ d\right]
}M_{k}^{-1}M_{f}\left( \exp ^{\left[ q-n\right] }r\right) }\leq
\frac{\rho _{h}^{\left( b,n\right) }\left( f\circ g\right) }{\lambda
_{k}^{\left( d,q\right) }\left( f\right)},
\end{equation*}
and\newline
$\noindent $\textbf{(ii)} If $q>m$, $a=c=p$, $\lambda _{f}\left( p,q\right)
>0$, $0<\lambda _{h}^{\left( b,q+n-m\right) }\left( f\circ g\right) \leq
\rho _{h}^{\left( b,q+n-m\right) }\left( f\circ g\right) $ $<\infty $ and $%
0<\lambda _{h}^{\left( d,q\right) }\left( f\right) \leq \rho _{k}^{\left(
d,q\right) }\left( f\right) <\infty$, then%
\begin{equation*}
\frac{\lambda _{h}^{\left( b,q+n-m\right) }\left( f\circ g\right) }{\rho
_{k}^{\left( d,q\right) }\left( f\right) }\leq \underset{r\rightarrow \infty
}{\lim \inf }\frac{\log ^{\left[ b\right] }M_{h}^{-1}M_{f\circ g}\left(
r\right) }{\log ^{\left[ d\right] }M_{k}^{-1}M_{f}\left( \exp ^{\left[ m-n%
\right] }r\right) }\leq \frac{\lambda _{h}^{\left( b,q+n-m\right) }\left(
f\circ g\right) }{\lambda _{k}^{\left( d,q\right) }\left( f\right) }~
\end{equation*}%
\begin{equation*}
~\ \ \ \ \ \ \ \ \ \ \ \ \ \ \ \ \ \ \ \ \ \ \leq \underset{r\rightarrow
\infty }{\lim \sup }\frac{\log ^{\left[ b\right] }M_{h}^{-1}M_{f\circ
g}\left( r\right) }{\log ^{\left[ d\right] }M_{k}^{-1}M_{f}\left( \exp ^{%
\left[ m-n\right] }r\right) }\leq \frac{\rho _{h}^{\left(
b,q+n-m\right) }\left( f\circ g\right) }{\lambda _{k}^{\left(
d,q\right) }\left( f\right)}.
\end{equation*}\end{theorem}

\begin{Proof}
Assume, that either ($q=m,$ $a=c=p,$ $q\geq n$) or ($q<m,$ $c=p,$ $a=p+m-q,$ $q\geq n$%
) hold and $\lambda _{f}\left( p,q\right) >0$. Then in view of Theorem \ref%
{t1}, the index-pair of $f\circ g$ is $\left( p,n\right)$ or $\left(
p+m-q,n\right)$, resp., and therefore by Definition \ref{d4}, $\rho
_{h}^{\left( b,n\right)}$ $\left( f\circ g\right)$ $\left( \lambda
_{h}^{\left( b,n\right) }\left( f\circ g\right),
\text{resp.}\right)$, and $\rho _{k}^{\left( d,q\right) }\left(
f\right) $ $\left(\lambda _{k}^{\left( d,q\right) }\left( f\right),
\text{resp.} \right) $ exist.
\newline Now from the
definition of $\rho _{k}^{\left( d,q\right) }\left( f\right) $ and
$\lambda _{h}^{\left( b,n\right) }\left( f\circ g\right)$,  for
arbitrary positive $\varepsilon$, and for all sufficiently large
values of $r$, we have
\begin{equation}
\log ^{\left[ b\right] }M_{h}^{-1}M_{f\circ g}\left( r\right) \geq
\left( \lambda _{h}^{\left( b,n\right) }\left( f\circ g\right)
-\varepsilon \right) \log ^{\left[ n\right] }r  \label{5.11}
\end{equation}%
and%
\begin{equation}
\log ^{\left[ d\right] }M_{k}^{-1}M_{f}\left( \exp ^{\left[
q-n\right] }r\right) \leq \left( \rho _{k}^{\left( d,q\right)
}\left( f\right) +\varepsilon \right) \log ^{\left[ n\right] }r.
\label{5.12}
\end{equation}%
Now from \eqref{5.11} and \eqref{5.12}, it follows for all
sufficiently large values of $r$, that
\begin{equation*}
\frac{\log ^{\left[ b\right] }M_{h}^{-1}M_{f\circ g}\left( r\right) }{\log ^{%
\left[ d\right] }M_{k}^{-1}M_{f}\left( \exp ^{\left[ q-n\right] }r\right) }%
\geq \frac{\left( \lambda _{h}^{\left( b,n\right) }\left( f\circ
g\right) -\varepsilon \right) \log ^{\left[ n\right] }r}{\left( \rho
_{k}^{\left( d,q\right) }\left( f\right) +\varepsilon \right) \log ^{\left[ n%
\right] }r}.
\end{equation*}%
Since $\varepsilon \left(\varepsilon >0\right) $ is arbitrary, we obtain that%
\begin{equation}
\underset{r\rightarrow \infty }{\lim \inf }\frac{\log ^{\left[ b\right]
}M_{h}^{-1}M_{f\circ g}\left( r\right) }{\log ^{\left[ d\right]
}M_{k}^{-1}M_{f}\left( \exp ^{\left[ q-n\right] }r\right) }\geq \frac{%
\lambda _{h}^{\left( b,n\right) }\left( f\circ g\right) }{\rho
_{k}^{\left( d,q\right) }\left( f\right)}.  \label{5.13}
\end{equation}%
For a sequence of values of $r$ tending to infinity we have
\begin{equation}
\log ^{\left[ b\right] }M_{h}^{-1}M_{f\circ g}\left( r\right) \leq
\left( \lambda _{h}^{\left( b,n\right) }\left( f\circ g\right)
+\varepsilon \right) \log ^{\left[ n\right] }r,  \label{5.14}
\end{equation}%
and for all sufficiently large values of $r$
\begin{equation}
\log ^{\left[ d\right] }M_{k}^{-1}M_{f}\left( \exp ^{\left[
q-n\right] }r\right) \geq \left( \lambda _{k}^{\left( d,q\right)
}\left( f\right) -\varepsilon \right) \log ^{\left[ n\right] }r.
\label{5.15}
\end{equation}%
Combining \eqref{5.14} and \eqref{5.15}, for a sequence of values of
$r$ tending to infinity, we get
\begin{equation*}
\frac{\log ^{\left[ b\right] }M_{h}^{-1}M_{f\circ g}\left( r\right) }{\log ^{%
\left[ d\right] }M_{k}^{-1}M_{f}\left( \exp ^{\left[ q-n\right] }r\right) }%
\leq \frac{\left( \lambda _{h}^{\left( b,n\right) }\left( f\circ
g\right) +\varepsilon \right) \log ^{\left[ n\right] }r}{\left(
\lambda _{k}^{\left( d,q\right) }\left( f\right) -\varepsilon
\right) \log ^{\left[ n\right] }r}.
\end{equation*}%
For arbitrary $\varepsilon \left(\varepsilon  >0\right)$, it follows
\begin{equation}
\underset{r\rightarrow \infty }{\lim \inf }\frac{\log ^{\left[ b\right]
}M_{h}^{-1}M_{f\circ g}\left( r\right) }{\log ^{\left[ d\right]
}M_{k}^{-1}M_{f}\left( \exp ^{\left[ q-n\right] }r\right) }\leq \frac{%
\lambda _{h}^{\left( b,n\right) }\left( f\circ g\right) }{\lambda
_{k}^{\left( d,q\right) }\left( f\right) }.  \label{5.16}
\end{equation}%
Also, for a sequence of values of $r$ tending to infinity, we obtain
\begin{equation}
\log ^{\left[ d\right] }M_{k}^{-1}M_{f}\left( \exp ^{\left[
q-n\right] }r\right) \leq \left( \lambda _{k}^{\left( d,q\right)
}\left( f\right) +\varepsilon \right) \log ^{\left[ n\right] }r.
\label{5.17}
\end{equation}%
Applying \eqref{5.11} and \eqref{5.17}, for a sequence of values of
$r$ tending to infinity, we get
\begin{equation*}
\frac{\log ^{\left[ b\right] }M_{h}^{-1}M_{f\circ g}\left( r\right) }{\log ^{%
\left[ d\right] }M_{k}^{-1}M_{f}\left( \exp ^{\left[ q-n\right] }r\right) }%
\geq \frac{\left( \lambda _{h}^{\left( b,n\right) }\left( f\circ
g\right) -\varepsilon \right) \log ^{\left[ n\right] }r}{\left(
\lambda _{k}^{\left( d,q\right) }\left( f\right) +\varepsilon
\right) \log ^{\left[ n\right] }r}.
\end{equation*}%
As $\varepsilon \left(\varepsilon >0\right)$ is arbitrary, we get from above that%
\begin{equation}
\underset{r\rightarrow \infty }{\lim \sup }\frac{\log ^{\left[ b\right]
}M_{h}^{-1}M_{f\circ g}\left( r\right) }{\log ^{\left[ d\right]
}M_{k}^{-1}M_{f}\left( \exp ^{\left[ q-n\right] }r\right) }\geq \frac{%
\lambda _{h}^{\left( b,n\right) }\left( f\circ g\right) }{\lambda
_{k}^{\left( d,q\right) }\left( f\right) }.  \label{5.18}
\end{equation}%
For all sufficiently large values of $r$ we obtain
\begin{equation}
\log T_{h}^{-1}T_{f\circ g}\left( r\right) \leq \left( \rho
_{h}^{\left(
b,n\right) }\left( f\circ g\right) +\varepsilon \right) \log ^{\left[ n%
\right] }r.  \label{5.19}
\end{equation}%
Combining now \eqref{5.15} and \eqref{5.19}, it follows for all
sufficiently large values of $r$
\begin{equation*}
\frac{\log ^{\left[ b\right] }M_{h}^{-1}M_{f\circ g}\left( r\right) }{\log ^{%
\left[ d\right] }M_{k}^{-1}M_{f}\left( \exp ^{\left[ q-n\right] }r\right) }%
\leq \frac{\left( \rho _{h}^{\left( b,n\right) }\left( f\circ
g\right) +\varepsilon \right) \log ^{\left[ n\right] }r}{\left(
\lambda _{k}^{\left( d,q\right) }\left( f\right) -\varepsilon
\right) \log ^{\left[ n\right] }r},
\end{equation*} and, therefore, for arbitrary $\varepsilon \left( >0\right)$, we obtain
\begin{equation}
\underset{r\rightarrow \infty }{\lim \sup }\frac{\log ^{\left[
b\right] }M_{h}^{-1}M_{f\circ g}\left( r\right) }{\log ^{\left[
d\right] }M_{k}^{-1}M_{f}\left( \exp ^{\left[ q-n\right] }r\right)
}\leq \frac{\rho _{h}\left( f\rho _{h}^{\left( b,n\right) }\left(
f\circ g\right) \circ g\right) }{\lambda _{k}^{\left( d,q\right)
}\left( f\right)}.  \label{5.20}
\end{equation}%
Thus the first part of the theorem follows from \eqref{5.13},
\eqref{5.16}, \eqref{5.18}, and \eqref{5.20}.\newline Similarly, one
can easily derive the second part of the theorem.
\end{Proof}{\hfill{$\Box$}{\medskip}}

Reasoning along the same line as in the proof of the Theorem
\ref{t3} we obtain:

\begin{theorem}
\label{t4} Let $f,g,h$ and $l$ be any four entire functions with
index-pairs $\left( p,q\right) $, $\left( m,n\right) ,$ $\left(
a,b\right) $ and $\left( x,y\right)$, resp., where $a,b,p,q,m,n,x,y$
are all positive integers such that $a\geq b,$ $p\geq q,m\geq n$ and
$x\geq y.$\newline
$\noindent $\textbf{(i)} If either $\left( q=m=x,\text{ }a=p\right) $ or $%
\left( q<m=x,\text{ }a=p+m-q\right) $ holds and $\lambda _{g}\left(
m,n\right) >0$, $0<\lambda _{h}^{\left( b,n\right) }\left( f\circ g\right)
\leq \rho _{h}^{\left( b,n\right) }\left( f\circ g\right) <\infty $, $%
0<\lambda _{l}^{\left( y,n\right) }\left( g\right) \leq \rho
_{l}^{\left( y,n\right) }\left( g\right) <\infty$, then
\begin{equation*}
\frac{\lambda _{h}^{\left( b,n\right) }\left( f\circ g\right) }{\rho
_{l}^{\left( y,n\right) }\left( g\right) }\leq \underset{r\rightarrow \infty
}{\lim \inf }\frac{\log ^{\left[ b\right] }M_{h}^{-1}M_{f\circ g}\left(
r\right) }{\log ^{\left[ y\right] }M_{l}^{-1}M_{g}\left( r\right) }\leq
\frac{\lambda _{h}^{\left( b,n\right) }\left( f\circ g\right) }{\lambda
_{l}^{\left( y,n\right) }\left( g\right) }~~\ \ \ \ \ \ \ \ \ \ \ \ \ \ \ \
\ \ \ \
\end{equation*}%
\begin{equation*}
~\ \ \ \ \ \ \ \ \ \ \ \ \ \ \ \ \ \ \ \ \ \ \ \ \ \ \ \ \ \ \ \ \ \
\ \ \ \leq \underset{r\rightarrow \infty }{\lim \sup }\frac{\log
^{\left[ b\right] }M_{h}^{-1}M_{f\circ g}\left( r\right) }{\log
^{\left[ y\right] }M_{l}^{-1}M_{g}\left( r\right) }\leq \frac{\rho
_{h}^{\left( b,n\right) }\left( f\circ g\right) }{\lambda
_{l}^{\left( y,n\right) }\left( g\right)},
\end{equation*}%
and\newline
$\noindent $\textbf{(ii)} If $q>m=x$, $a=p$, $\lambda _{g}\left( m,n\right)
>0$, $0<\lambda _{h}^{\left( b,q+n-m\right) }\left( f\circ g\right) \leq
\rho _{h}^{\left( b,q+n-m\right) }\left( f\circ g\right) <\infty $, $%
0<\lambda _{l}^{\left( y,n\right) }\left( g\right) \leq \rho
_{l}^{\left( y,n\right) }\left( g\right) <\infty$, then
\begin{equation*}
\frac{\lambda _{h}^{\left( b,q+n-m\right) }\left( f\circ g\right) }{\rho
_{l}^{\left( y,n\right) }\left( g\right) }\leq \underset{r\rightarrow \infty
}{\lim \inf }\frac{\log ^{\left[ b\right] }M_{h}^{-1}M_{f\circ g}\left( \exp
^{\left[ q-m\right] }r\right) }{\log ^{\left[ y\right] }M_{l}^{-1}M_{g}%
\left( r\right) }\leq \frac{\lambda _{h}^{\left( b,q+n-m\right) }\left(
f\circ g\right) }{\lambda _{l}^{\left( y,n\right) }\left( g\right) }
\end{equation*}%
\begin{equation*}
~\ \ \ \ \ \ \ \ \ \ \ \ \ \ \ \ \ \ \ \ \ \ \leq
\underset{r\rightarrow \infty }{\lim \sup }\frac{\log ^{\left[
b\right] }M_{h}^{-1}M_{f\circ g}\left( \exp ^{\left[ q-m\right]
}r\right) }{\log ^{\left[ y\right] }M_{l}^{-1}M_{g}\left( r\right)
}\leq \frac{\rho _{h}^{\left( b,q+n-m\right) }\left( f\circ g\right)
}{\lambda _{l}^{\left( y,n\right) }\left( g\right)}.
\end{equation*}
\end{theorem}

\begin{theorem}
\label{t5} Let $f,g,h$ and $k$ be any four entire functions with
index-pairs $\left( p,q\right) $, $\left( m,n\right) ,$ $\left(
a,b\right) $ and $\left( c,d\right)$, resp., where $a,b,c,d,p,q,m,n$
are all positive integers with $a\geq b$, $c\geq d,$ $p\geq q$ and
$m\geq n.$\newline
$\noindent $\textbf{(i)} If either ($q=m,$ $a=c=p,$ $q\geq n$) or ($q<m,$ $%
c=p,$ $a=p+m-q,$ $q\geq n$) holds and $\lambda _{f}\left( p,q\right) >0$, $%
0<\rho _{h}^{\left( b,n\right) }\left( f\circ g\right) <\infty $,
$0<\rho _{k}^{\left( d,q\right) }\left( f\right) <\infty$, then
\begin{equation*}
\underset{r\rightarrow \infty }{\lim \inf }\frac{\log ^{\left[ b\right]
}M_{h}^{-1}M_{f\circ g}\left( r\right) }{\log ^{\left[ d\right]
}M_{k}^{-1}M_{f}\left( \exp ^{\left[ q-n\right] }r\right) }\leq \frac{\rho
_{h}^{\left( b,n\right) }\left( f\circ g\right) }{\rho _{k}^{\left(
d,q\right) }\left( f\right) }\leq \underset{r\rightarrow \infty }{\lim \sup }%
\frac{\log ^{\left[ b\right] }M_{h}^{-1}M_{f\circ g}\left( r\right) }{\log ^{%
\left[ d\right] }M_{k}^{-1}M_{f}\left( \exp ^{\left[ q-n\right]
}r\right)},
\end{equation*}
and\newline
$\noindent $\textbf{(ii)} If $q>m$, $a=c=p$, $\lambda _{f}\left( p,q\right)
>0$, $0<\rho _{h}^{\left( b,q+n-m\right) }\left( f\circ g\right) $ $<\infty $
and $0<\rho _{k}^{\left( d,q\right) }\left( f\right) <\infty$, then
\begin{eqnarray*}
\underset{r\rightarrow \infty }{\lim \inf }\frac{\log ^{\left[ b\right]
}M_{h}^{-1}M_{f\circ g}\left( r\right) }{\log ^{\left[ d\right]
}M_{k}^{-1}M_{f}\left( \exp ^{\left[ m-n\right] }r\right) } &\leq &\frac{%
\rho _{h}^{\left( b,q+n-m\right) }\left( f\circ g\right) }{\rho _{k}^{\left(
d,q\right) }\left( f\right) } \\
&\leq &\underset{r\rightarrow \infty }{\lim \sup }\frac{\log ^{\left[ b%
\right] }M_{h}^{-1}M_{f\circ g}\left( r\right) }{\log ^{\left[
d\right] }M_{k}^{-1}M_{f}\left( \exp ^{\left[ m-n\right]}r\right)}.
\end{eqnarray*}
\end{theorem}

\begin{Proof}
Let either ($q=m,$ $a=c=p,$ $q\geq n$) or ($q<m,$ $c=p,$ $a=p+m-q,$ $q\geq n$%
) hold, and also let $\lambda _{f}\left( p,q\right) >0$. In view of
Theorem \ref{t1}, the index-pair of $f\circ g\ $ is $\ \left( p,n\right) $ or $%
\left( p+m-q,n\right)$, resp. Hence by Definition \ref{d4}, $\rho
_{h}^{\left( b,n\right) }\left( f\circ g\right) $ and $\rho
_{k}^{\left( d,q\right) }\left( f\right)$ exist, and from the
definition of $\rho _{k}^{\left( d,q\right) }\left( f\right)$, for a
sequence of values of $r$ tending to infinity, we get
\begin{eqnarray}
\log ^{\left[ d\right] }M_{k}^{-1}M_{f}\left( \exp ^{\left[
q-n\right] }r\right) &\geq &\left( \rho _{k}^{\left( d,q\right)
}\left( f\right)
-\varepsilon \right) \log ^{\left[ n\right] }r  \notag  \label{3.1} \\
i.e.,~\log T_{P[h]}^{-1}T_{P[f]}\left( r\right) &\geq &\left( \rho
_{k}^{\left( d,q\right) }\left( f\right) -\varepsilon \right) \log ^{\left[ n%
\right] }r.  \label{5.21}
\end{eqnarray}%
Now from \eqref{5.19} and \eqref{5.21}, for a sequence of values of
$r$ tending to infinity, it follows  that%
\begin{equation*}
\frac{\log ^{\left[ b\right] }M_{h}^{-1}M_{f\circ g}\left( r\right) }{\log ^{%
\left[ d\right] }M_{k}^{-1}M_{f}\left( \exp ^{\left[ q-n\right] }r\right) }%
\leq \frac{\left( \rho _{h}^{\left( b,n\right) }\left( f\circ
g\right) +\varepsilon \right) \log ^{\left[ n\right] }r}{\left( \rho
_{k}^{\left( d,q\right) }\left( f\right) -\varepsilon \right) \log
^{\left[ n\right] }r}.
\end{equation*}%
As $\varepsilon \left( >0\right) $ is arbitrary, we obtain that%
\begin{equation}
\underset{r\rightarrow \infty }{\lim \inf }\frac{\log ^{\left[
b\right] }M_{h}^{-1}M_{f\circ g}\left( r\right) }{\log ^{\left[
d\right] }M_{k}^{-1}M_{f}\left( \exp ^{\left[ q-n\right] }r\right)
}\leq \frac{\rho _{h}^{\left( b,n\right) }\left( f\circ g\right)
}{\rho _{k}^{\left( d,q\right) }\left( f\right) }.  \label{5.22}
\end{equation}%
Again,  for a sequence of values of $r$ tending to infinity, we
obtain
\begin{equation}
\log ^{\left[ b\right] }M_{h}^{-1}M_{f\circ g}\left( r\right) \geq
\left( \rho _{h}^{\left( b,n\right) }\left( f\circ g\right)
-\varepsilon \right) \log ^{\left[ n\right] }r.  \label{5.23}
\end{equation}%
Combining \eqref{5.12} and  \eqref{5.23}, for a sequence of values
of $r$ tending to infinity, we get
\begin{equation*}
\frac{\log ^{\left[ b\right] }M_{h}^{-1}M_{f\circ g}\left( r\right) }{\log ^{%
\left[ d\right] }M_{k}^{-1}M_{f}\left( \exp ^{\left[ q-n\right] }r\right) }%
\geq \frac{\left( \rho _{h}^{\left( b,n\right) }\left( f\circ
g\right) -\varepsilon \right) \log ^{\left[ n\right] }r}{\left( \rho
_{k}^{\left( d,q\right) }\left( f\right) +\varepsilon \right) \log
^{\left[ n\right] }r}.
\end{equation*}%
For arbitrarily chosen  $\varepsilon \left( >0\right)$, it follows
from the above
\begin{equation}
\underset{r\rightarrow \infty }{\lim \sup }\frac{\log ^{\left[ b\right]
}M_{h}^{-1}M_{f\circ g}\left( r\right) }{\log ^{\left[ d\right]
}M_{k}^{-1}M_{f}\left( \exp ^{\left[ q-n\right] }r\right) }\geq \frac{%
\rho _{h}^{\left( b,n\right) }\left( f\circ g\right) }{\rho
_{k}^{\left( d,q\right) }\left( f\right)}.  \label{5.24}
\end{equation}%
Thus the first part of the theorem follows from \eqref{5.22} and
\eqref{5.24}. \newline Analogously, the second part of the proof of
the theorem can be derived.
\end{Proof}{\hfill{$\Box$}{\medskip}}

The proof of the following theorem can be carried out as of the
Theorem \ref {t5}, therefore we omit the details.

\begin{theorem}
\label{t6} Let $f,g,h$ and $l$ be any four entire functions with
index-pairs $\left( p,q\right) $, $\left( m,n\right) ,$ $\left(
a,b\right) $ and $\left( x,y\right)$, resp., where $a,b,p,q,m,n,x,y$
are all positive integers such that $a\geq b,$ $p\geq q,m\geq n$ and
$x\geq y.$\newline
$\noindent $\textbf{(i)} If either $\left( q=m=x,\text{ }a=p\right) $ or $%
\left( q<m=x,\text{ }a=p+m-q\right) $ holds and $\lambda _{g}\left(
m,n\right) >0$, $0<\rho _{h}^{\left( b,n\right) }\left( f\circ
g\right) <\infty $, $0<\rho _{l}^{\left( y,n\right) }\left( g\right)
<\infty$, then
\begin{equation*}
\underset{r\rightarrow \infty }{\lim \inf }\frac{\log ^{\left[ b\right]
}M_{h}^{-1}M_{f\circ g}\left( r\right) }{\log ^{\left[ y\right]
}M_{l}^{-1}M_{g}\left( r\right) }\leq \frac{\rho _{h}^{\left( b,n\right)
}\left( f\circ g\right) }{\rho _{l}^{\left( y,n\right) }\left( g\right) }%
\leq \underset{r\rightarrow \infty }{\lim \sup }\frac{\log ^{\left[
b\right] }M_{h}^{-1}M_{f\circ g}\left( r\right) }{\log ^{\left[
y\right] }M_{l}^{-1}M_{g}\left( r\right)},
\end{equation*}%
and\newline
$\noindent $\textbf{(ii)} If $q>m=x$, $a=p$, $\lambda _{g}\left( m,n\right)
>0$, $0<\rho _{h}^{\left( b,q+n-m\right) }\left( f\circ g\right) <\infty $, $%
0<\rho _{l}^{\left( y,n\right) }\left( g\right) <\infty$, then%
\begin{eqnarray*}
\underset{r\rightarrow \infty }{\lim \inf }\frac{\log ^{\left[ b\right]
}M_{h}^{-1}M_{f\circ g}\left( \exp ^{\left[ q-m\right] }r\right) }{\log ^{%
\left[ y\right] }M_{l}^{-1}M_{g}\left( r\right) } &\leq &\frac{\rho
_{h}^{\left( b,q+n-m\right) }\left( f\circ g\right) }{\rho _{l}^{\left(
y,n\right) }\left( g\right) } \\
&\leq &\underset{r\rightarrow \infty }{\lim \sup }\frac{\log ^{\left[ b%
\right] }M_{h}^{-1}M_{f\circ g}\left( \exp ^{\left[ q-m\right] }r\right) }{%
\log ^{\left[ y\right] }M_{l}^{-1}M_{g}\left( r\right)}.
\end{eqnarray*}
\end{theorem}

The following theorem is a natural consequence of Theorem \ref{t3}
and Theorem \ref{t5}:

\begin{theorem}
\label{t7} Let $f,g,h$ and $k$ be any four entire functions with
index-pairs $\left( p,q\right) $, $\left( m,n\right) ,$ $\left(
a,b\right)$ (and $\left( c,d\right)$, resp.), where
$a,b,c,d,p,q,m,n$ are all positive integers such that $a\geq b$,
$c\geq d,$ $p\geq q$ and $m\geq n.$\newline $\noindent $\textbf{(i)}
If either $\left( q=m,\text{ }a=c=p,\text{ }q\geq n\right) $ or
$\left( q<m,\text{ }c=p,\text{ }a=p+m-q,\text{ }q\geq n\right) $
holds and $\lambda _{f}\left( p,q\right) >0$, $0<\lambda
_{h}^{\left( b,n\right) }\left( f\circ g\right) \leq \rho
_{h}^{\left( b,n\right) }\left( f\circ g\right) <\infty $,
$0<\lambda _{h}^{\left( d,q\right) }\left(
f\right) \leq \rho _{k}^{\left( d,q\right) }\left( f\right) <\infty$, then%
\begin{equation*}
\underset{r\rightarrow \infty }{\lim \inf }\frac{\log ^{\left[ b\right]
}M_{h}^{-1}M_{f\circ g}\left( r\right) }{\log ^{\left[ d\right]
}M_{k}^{-1}M_{f}\left( \exp ^{\left[ q-n\right] }r\right) }\leq \min \left\{
\frac{\lambda _{h}^{\left( b,n\right) }\left( f\circ g\right) }{\lambda
_{k}^{\left( d,q\right) }\left( f\right) },\frac{\rho _{h}^{\left(
b,n\right) }\left( f\circ g\right) }{\rho _{k}^{\left( d,q\right) }\left(
f\right) }\right\} \leq ~\ \ \
\end{equation*}%
\begin{equation*}
~\ \ \ \ \ \ \ \ \ \max \left\{ \frac{\lambda _{h}^{\left( b,n\right)
}\left( f\circ g\right) }{\lambda _{k}^{\left( d,q\right) }\left( f\right) },%
\frac{\rho _{h}^{\left( b,n\right) }\left( f\circ g\right) }{\rho
_{k}^{\left( d,q\right) }\left( f\right) }\right\} \leq \underset{%
r\rightarrow \infty }{\lim \sup }\frac{\log ^{\left[ b\right]
}M_{h}^{-1}M_{f\circ g}\left( r\right) }{\log ^{\left[ d\right]
}M_{k}^{-1}M_{f}\left( \exp ^{\left[ q-n\right] }r\right)},
\end{equation*}%
and\newline
$\noindent $\textbf{(ii)} If $q>m$, $a=c=p$, $\lambda _{f}\left( p,q\right)
>0$, $0<\lambda _{h}^{\left( b,q+n-m\right) }\left( f\circ g\right) \leq
\rho _{h}^{\left( b,q+n-m\right) }\left( f\circ g\right) $ $<\infty $ and $%
0<\lambda _{h}^{\left( d,q\right) }\left( f\right) \leq \rho _{k}^{\left(
d,q\right) }\left( f\right) <\infty $, then%
\begin{equation*}
\underset{r\rightarrow \infty }{\lim \inf }\frac{\log ^{\left[ b\right]
}M_{h}^{-1}M_{f\circ g}\left( r\right) }{\log ^{\left[ d\right]
}M_{k}^{-1}M_{f}\left( \exp ^{\left[ m-n\right] }r\right) }\leq \min \left\{
\frac{\lambda _{h}^{\left( b,q+n-m\right) }\left( f\circ g\right) }{\lambda
_{k}^{\left( d,q\right) }\left( f\right) },\frac{\rho _{h}^{\left(
b,q+n-m\right) }\left( f\circ g\right) }{\rho _{k}^{\left( d,q\right)
}\left( f\right) }\right\} \leq
\end{equation*}%
\begin{equation*}
\max \left\{ \frac{\lambda _{h}^{\left( b,q+n-m\right) }\left( f\circ
g\right) }{\lambda _{k}^{\left( d,q\right) }\left( f\right) },\frac{\rho
_{h}^{\left( b,q+n-m\right) }\left( f\circ g\right) }{\rho _{k}^{\left(
d,q\right) }\left( f\right) }\right\} \leq \underset{r\rightarrow \infty }{%
\lim \sup }\frac{\log ^{\left[ b\right] }M_{h}^{-1}M_{f\circ g}\left(
r\right) }{\log ^{\left[ d\right] }M_{k}^{-1}M_{f}\left( \exp ^{\left[ m-n%
\right] }r\right) }.
\end{equation*}
\end{theorem}

We omit the proof because of the similarity to the previous ones.

Analogously one may formulate the following theorem without its
proof.

\begin{theorem}
\label{t8} Let $f,g,h$ and $l$ be any four entire functions with
index-pairs $\left( p,q\right) $, $\left( m,n\right) ,$ $\left(
a,b\right)$ (and $\left( x,y\right)$, resp.), where
$a,b,p,q,m,n,x,y$ are all positive integers such that $a\geq b,$
$p\geq q,m\geq n$ and $x\geq y.$\newline
$\noindent $\textbf{(i)} If either $\left( q=m=x,\text{ }a=p\right) $ or $%
\left( q<m=x,\text{ }a=p+m-q\right) $ holds and $\lambda _{g}\left(
m,n\right) >0$, $0<\lambda _{h}^{\left( b,n\right) }\left( f\circ g\right)
\leq \rho _{h}^{\left( b,n\right) }\left( f\circ g\right) <\infty $, $%
0<\lambda _{l}^{\left( y,n\right) }\left( g\right) \leq \rho _{l}^{\left(
y,n\right) }\left( g\right) <\infty $, then%
\begin{equation*}
\underset{r\rightarrow \infty }{\lim \inf }\frac{\log ^{\left[ b\right]
}M_{h}^{-1}M_{f\circ g}\left( r\right) }{\log ^{\left[ y\right]
}M_{l}^{-1}M_{g}\left( r\right) }\leq \min \left\{ \frac{\lambda
_{h}^{\left( b,n\right) }\left( f\circ g\right) }{\lambda _{l}^{\left(
y,n\right) }\left( g\right) },\frac{\rho _{h}^{\left( b,n\right) }\left(
f\circ g\right) }{\rho _{l}^{\left( y,n\right) }\left( g\right) }\right\}
\leq ~\ \ \ \ \ \ \ \ \ \
\end{equation*}%
\begin{equation*}
~~\ \ \ \ \ \ \ \ \ \ \ \ \max \left\{ \frac{\lambda _{h}^{\left( b,n\right)
}\left( f\circ g\right) }{\lambda _{l}^{\left( y,n\right) }\left( g\right) },%
\frac{\rho _{h}^{\left( b,n\right) }\left( f\circ g\right) }{\rho
_{l}^{\left( y,n\right) }\left( g\right) }\right\} \leq \underset{%
r\rightarrow \infty }{\lim \sup }\frac{\log ^{\left[ b\right]
}M_{h}^{-1}M_{f\circ g}\left( r\right) }{\log ^{\left[ y\right]
}M_{l}^{-1}M_{g}\left( r\right) },
\end{equation*}%
and\newline
$\noindent $\textbf{(ii)} If $q>m=x$, $a=p$, $\lambda _{g}\left( m,n\right)
>0$, $0<\lambda _{h}^{\left( b,q+n-m\right) }\left( f\circ g\right) \leq
\rho _{h}^{\left( b,q+n-m\right) }\left( f\circ g\right) <\infty $, $%
0<\lambda _{l}^{\left( y,n\right) }\left( g\right) \leq \rho _{l}^{\left(
y,n\right) }\left( g\right) <\infty $, then%
\begin{equation*}
\underset{r\rightarrow \infty }{\lim \inf }\frac{\log ^{\left[ b\right]
}M_{h}^{-1}M_{f\circ g}\left( \exp ^{\left[ q-m\right] }r\right) }{\log ^{%
\left[ y\right] }M_{l}^{-1}M_{g}\left( r\right) }\leq \min \left\{ \frac{%
\lambda _{h}^{\left( b,q+n-m\right) }\left( f\circ g\right) }{\lambda
_{l}^{\left( y,n\right) }\left( g\right) },\frac{\rho _{h}^{\left(
b,q+n-m\right) }\left( f\circ g\right) }{\rho _{l}^{\left( y,n\right)
}\left( g\right) }\right\} \leq
\end{equation*}%
\begin{equation*}
\max \left\{ \frac{\lambda _{h}^{\left( b,q+n-m\right) }\left( f\circ
g\right) }{\lambda _{l}^{\left( y,n\right) }\left( g\right) },\frac{\rho
_{h}^{\left( b,q+n-m\right) }\left( f\circ g\right) }{\rho _{l}^{\left(
y,n\right) }\left( g\right) }\right\} \leq \underset{r\rightarrow \infty }{%
\lim \sup }\frac{\log ^{\left[ b\right] }M_{h}^{-1}M_{f\circ g}\left( \exp ^{%
\left[ q-m\right] }r\right) }{\log ^{\left[ y\right]
}M_{l}^{-1}M_{g}\left( r\right) }.
\end{equation*}
\end{theorem}\bigskip

\centerline{ \textbf{Acknowledgements}}

\noindent This work was partially supported by the Centre for
Innovation and Transfer of Natural Sciences and Engineering
Knowledge, Faculty of Mathematics and Natural Sciences, University
of Rzeszow.

\end{document}